\documentclass[11pt]{amsart}
mm \hoffset -16 true mm

\begin{document}
\title{the automatic additivity of  $\xi-$Lie derivations on von Neumann algebras }
\author{Zhaofang Bai}
\address[Zhaofang Bai]{School of Mathematical Sciences, Xiamen University, Xiamen,
361005, P. R. China.} \email[Zhaofang
Bai]{zhaofangbai@yahoo.com.cn}
\author{Shuanping Du$^*$}
\address[Shuanping Du]{School of Mathematical Sciences, Xiamen University, Xiamen, 361005,
 P. R. China.} \email[Shuanping
Du]{shuanpingdu@yahoo.com}
\author{Yu Guo}
\address[Yu Guo]{Department of Mathematics, Shanxi Datong University, Datong, 037009, P. R. China}
\email[Yu Guo]{guoyu3@yahoo.com.cn}
\thanks{This work was supported partially by National Natural Science
Foundation of China (11071201, 11001230), and the Fundamental
Research Funds for the Central Universities (2010121001).}
\thanks{This paper is in final form and no version of it will be submitted for
publication elsewhere.}
\thanks{$^*$ Correspondence author }
\begin{abstract}
Let ${\mathcal M}$  be a von Neumann algebra with no central
summands of type I$_1$. It is shown that every  nonlinear
$\xi-$Lie derivation ($\xi\neq 1$) on $\mathcal M$ is an additive
derivation. \end{abstract} \maketitle

\section{Introduction and main results}
Let $\mathcal A$ be an associate ring (or an algebra over a field
$\mathbb F$). Then ${\mathcal A}$ is a Lie ring (Lie algebra)
under the product $[x,y]=xy-yx$, i.e., the commutator of $x$ and
$y$. Recall that an additive (linear) map $\delta :{\mathcal
A}\rightarrow {\mathcal A}$ is called an additive (linear)
derivation if $\delta(xy)=\delta(x)y+x\delta(y)$ for all $x,y\in
{\mathcal A}$. Derivations are very important maps both in theory
and in applications, and have been studied intensively (see
\cite{Ch,Sa,Se1,Se2} and the references therein). More generally,
an additive (linear) map $L$ from ${\mathcal A}$ into itself is
called an additive (linear) Lie derivation if
$L([x,y])=[L(x),y]+[x,L(y)]$ for all $x,y\in {\mathcal A}$. The
questions of characterizing Lie derivations and revealing the
relationship between Lie derivations and derivations have received
many mathematicians' attention recently (see
\cite{Br2,John,MaVi,Mie3}). Very roughly speaking, additive
(linear) Lie derivations in the context prime rings (operator
algebras) can be decomposed as $\sigma +\tau$, where $\sigma$ is
an additive (linear) derivation and $\tau$ is an additive (linear)
map sending commutators into zero. Similarly, associated with the
Jordan product $xy+yx$. we have the conception of Jordan
derivation which is also studied intensively (see
\cite{Br3,Br4,John} and the references therein).

Note that an important relation associated with the Lie product is
the commutativity. Two elements $x,y$ in an algebra ${\mathcal A}$
are commutative if $xy=yx$, that is, their Lie product is zero.
More generally, if $\xi$ is a scalar and if $xy=\xi yx$, we say
that $x$ commutes with $y$ up to a factor $\xi$. The notion of
commutativity up to a factor for pairs of operators is also
important and has been studied in the context of operator algebras
and quantum groups (Refs. \cite{BrBuPe,Ka}). Motivated by this,
the authors introduce a binary operation $[x,y]_{\xi}=xy-\xi yx$,
called $\xi$-Lie product of $x,y$ (Ref. \cite{QiHo1}). This
product is found playing a more and more important role in some
research topics, and its study has recently attracted many
authors¡¯ attention (for example, see \cite{QiHo1,QiHo2}). Then it
is natural to introduce the concept of $\xi$-Lie derivation. An
additive (linear) map $L$ from ${\mathcal A}$  into itself is
called a $\xi-$Lie derivation if
$L([x,y]_{\xi})=[L(x),y]_{\xi}+[x,L(y)]_{\xi}$ for all $x,y\in
{\mathcal A}$. This concept unifies several well-known notions. It
is clear that a $\xi-$Lie derivation is a derivation if $\xi=0$;
is a Lie derivation if $\xi=1$; is a Jordan derivation if
$\xi=-1$. In \cite{QiHo2}, Qi and Hou characterized the additive
$\xi-$Lie derivation on nest algebras.

Let $\Phi:{\mathcal A}\rightarrow {\mathcal A}$ be a map (without
the additivity or linearity assumption). We say that $\Phi$ is a
nonlinear $\xi-$Lie derivation if
$\Phi([x,y]_{\xi})=[\Phi(x),y]_{\xi}+[x,\Phi(y)]_{\xi}$ for all
$x,y\in {\mathcal A}$. Recently, Yu and Zhang \cite{YuZh}
described nonlinear Lie derivation on triangular algebras. The aim
of this note is to investigate  nonlinear $\xi-$Lie derivations on
von Neumann algebras ($\xi\neq 1$) and to reveal the relationship
between such nonlinear $\xi-$Lie derivations and additive
derivations. Due to vital importance of derivations, we firstly
investigate nonlinear derivations. To our surprising,  nonlinear
derivations are automatically additive.  Our main results read as
follows.

\textbf{Theorem 1.1.} {\it Let ${\mathcal M}$ be a von Neumann
algebra with no central summands of type I$_1$.  If
$\Phi:{\mathcal M}\rightarrow{\mathcal M}$ is a nonlinear
derivation, then $\Phi$ is an additive derivation.}

The following result reveals the relationship between general
nonlinear $\xi-$Lie derivations and additive derivations.

\textbf{Theorem 1.2.} {\it Let ${\mathcal M}$ be a von Neumann
algebra with no central summands of type I$_1$.  If $\xi$ is a
scalar not equal $0,1$ and $\Phi:{\mathcal M}\rightarrow{\mathcal
M}$ is a nonlinear $\xi-$Lie derivation, then $\Phi$ is an
additive derivation and $\Phi(\xi T)=\xi\Phi(T)$ for all $T\in
{\mathcal M}$.}

It is worth mentioning that, as it turns out from Theorems 1.1 and
Theorem 1.2, the additive structure and $\xi-$Lie multiplicative
structure of von Neumann algebra with no central summands of type
I$_1$ are very closely related to each other. We remark that the
question when a multiplicative map is necessary additive is
important in quantum mechanics and mathematics, and was discussed
for associative rings in the purely algebraic setting (\cite{Mar},
for a recent systematic account, see \cite{BD}). In recent years,
there is a growing interest in studying the automatic additivity
of maps determined by the action on the product (see
\cite{AH,BD,Lu1,QH,Wang} and the references therein). We also
remark that if $\xi=1$, then $\xi-$Lie derivation is in fact a Lie
derivation, while Lie derivation is not necessary additive. For
example, let $\sigma$ is an additive derivation of ${\mathcal M}$
and $\tau$ is a mapping of ${\mathcal M}$ into its center
${\mathcal Z}_{\mathcal M}$ which maps commutators into zero. Then
$\sigma +\tau$ is a Lie derivation and such Lie derivation is not
additive in general.

\section{Notations and Preliminaries}

Before embarking on the proof of our main results, we need some
notations and preliminaries about von Neumann algebras. A von
Neumann algebra $\mathcal M$ is a weakly closed, self-adjoint
algebra of operators on a Hilbert space $H$ containing the
identity operator I. The set ${\mathcal Z}_{\mathcal
M}=\{S\in{\mathcal M}\mid ST=TS \text{ for all } T\in{\mathcal
M}\}$ is called the center of ${\mathcal M}$. For $A\in {\mathcal
M}$, the central carrier of $A$, denoted by $\overline{A}$, is the
intersection of all central projections $P$ such that $PA=A$. It
is well known that the central carrier of $A$ is the projection
with the range $[{\mathcal M}A(H)]$, the closed linear span of
$\{MA(x)\mid M\in {\mathcal M}, x\in H\}$. For each self-adjoint
operator $A\in {\mathcal M}$, we define the central core of $A$,
denoted by $\underline{A}$, to be $\sup\{S\in {\mathcal
Z}_{\mathcal M}\mid S=S^*, S\leq A\}$. Clearly, one has
$A-\underline{A}\geq 0$. Further if $S\in {\mathcal Z}_{\mathcal
M}$ and $A-\underline{A}\geq S\geq 0$ then $S=0$. If $P$ is a
projection it is clear that $\underline{P}$ is the largest central
projection $\leq P$. We call a projection core-free if
$\underline{P}=0$. It is easy to see that $\underline{P}=0$ if and
only if $\overline{I-P}=I$, here $\overline{I-P}$ denotes the
central carrier of $I-P$. We use \cite{Kad} as a general reference
for the theory of von Neumann algebras.

In the following, there are several fundamental  properties of von
Neumann algebras  from \cite{Br,Mie2} which will be used
frequently. For convenience, we list them in a lemma.

\textbf{Lemma 2.1.} Let $\mathcal M$ be a von Neumann algebra.

(i) (\cite[Lemma 4]{Mie2}) If  ${\mathcal M}$ has no summands of
type I$_1$, then each nonzero central projection of ${\mathcal M}$
is the central carrier of a core-free projection of ${\mathcal
M}$;

(ii) (\cite[Lemma 2.6]{Br}) If  ${\mathcal M}$ has no summands of
type I$_1$, then $\mathcal M$ equals the ideal of $\mathcal M$
generated by all commutators in $\mathcal M$.

By Lemma 2.1(i), one can find a non-trivial core-free projection
with central carrier $I$, denoted  by $P_1$. Throughout this
paper, $P_1$ is fixed. Write $P_2=I-P_1$.  By the definition of
central core and central carrier, $P_2$ is also core-free and
$\overline{P_2}=I$. According to the two-side Pierce decomposition
of $\mathcal M$ relative $P_1$, denote ${\mathcal
M}_{ij}=P_i{\mathcal M}P_j$, $i,j=1,2$, then we may write
${\mathcal M}={\mathcal M}_{11}+{\mathcal M}_{12}+{\mathcal
M}_{21}+ {\mathcal M}_{22}$. In all that follows, when we write
$T_{ij}$, $S_{ij}$, $M_{ij}$, it indicates that they are contained
in ${\mathcal M}_{ij}$. A conclusion which is used frequently is
$TM_{ij}=0$ for every $M_{ij}\in{\mathcal M}_{ij}$ implies that
$TP_i=0$. Indeed $TP_iMP_j=0$ for all $M\in {\mathcal M}$ together
with $\overline{P_j}=I$    gives  $TP_i=0$. Similarly, if
$M_{ij}T=0$ for every $M_{ij}\in{\mathcal M}_{ij}$, then
$T^*M_{ij}^*=0$ and so $P_jT=0$. If $Z\in {\mathcal Z}_{\mathcal
M}$ and $ZP_i=0$, then $ZMP_i=0$ for all $M\in {\mathcal M}$ which
implies $Z=0$.

The next lemma is technical which plays an important role in the
proof of Theorem 1.2.

\textbf{Lemma 2.2.} Let $T\in {\mathcal M}$, $\xi\neq  0, 1$. Then
$T\in {\mathcal M}_{ij}+(\xi P_i+P_j){\mathcal Z}_{\mathcal M}$
($1\leq i\neq j\leq 2$) if and only if $[T,M_{ij}]_{\xi}=0$ for
every $M_{ij}\in {\mathcal M}_{ij}$;

\textbf{Proof.} The necessity is clear. Conversely, assume
$[T,M_{ij}]_{\xi}=0$ for every $M_{ij}\in {\mathcal M}_{ij}$.
Write $T=\sum_{i,j=1}^{2}T_{ij}$. It follows that $T_{ii}M_{ij}+
T_{ji}M_{ij}=\xi(M_{ij}T_{jj}+M_{ij}T_{ji})$. Thus
$$T_{ii}M_{ij}=\xi M_{ij}T_{jj} \eqno{(1)}$$ and $T_{ji}M_{ij}=0$.
Noting that  $\overline{P_j}=I$, we obtain $$T_{ji}=0.$$
 For every $M_{ii}\in {\mathcal M}_{ii}$,  $M_{jj}\in {\mathcal M}_{jj}$, $M_{ii}M_{ij}, M_{ij}M_{jj} \in {\mathcal M}_{ij}$
 and so $TM_{ii}M_{ij}=\xi M_{ii}M_{ij}T$ and $TM_{ij}M_{jj}=\xi M_{ij}M_{jj}T$. From $[T,M_{ij}]_{\xi}=0$, it follows that
 $TM_{ii}M_{ij}= M_{ii}TM_{ij}$, that is $(TM_{ii}-M_{ii}T)M_{ij}=0$. Using $\overline{P_j}=I$ again,
 we have $T_{ii}M_{ii}-M_{ii}T_{ii}=0$, i.e., $T_{ii}\in {\mathcal Z}_{P_i{\mathcal M}P_i}$.
 Thus $$T_{ii}=Z_iP_i$$ for some central element $Z_i\in {\mathcal Z}_{\mathcal M}$.
 Similarly, combining $TM_{ij}M_{jj}=\xi M_{ij}M_{jj}T$ and $[T,M_{ij}]_{\xi}=0$, we can obtain $$T_{jj}=Z_jP_j$$ for some central element
 $Z_j\in {\mathcal Z}_{\mathcal M}$.
 Now equation (1) implies that $(Z_i-\xi Z_j)M_{ij}=0$.
 From $\overline{P_j}=I$ and $M_{ij}$ is arbitrary, it follows that $(Z_i-\xi Z_j)P_i=0$.
  Since $Z_i-\xi Z_j\in {\mathcal Z}_{\mathcal M}$, $M(Z_i-\xi Z_j)P_i=(Z_i-\xi Z_j)MP_i=0$  for all $M\in \mathcal M$. By $\overline{P_i}=I$,
  it follows that $Z_i=\xi Z_j$. So  $T=T_{ij}+(\xi P_i+ P_j)Z_j\in {\mathcal M}_{ij}+(\xi P_i+ P_j){\mathcal Z}_{\mathcal M}$.

\section{Proofs of main results}

In the following, we are firstly aimed to prove Theorem 1.1.

\textbf{Proof of Theorem 1.1.} In what follows, $\Phi: {\mathcal
M}\rightarrow {\mathcal M}$ is a nonlinear derivation. We will
prove that $\Phi$ is additive, that is, for all $T,S\in{\mathcal
M}$, $\Phi(T+S)=\Phi(T)+\Phi(S)$. It is clear that
$\Phi(0)=\Phi(0)0+0\Phi(0)=0$. Note that
$\Phi(P_1P_2)=\Phi(P_1)P_2+P_1\Phi(P_2)=0$, multiplying by $P_2$
from the both sides of this equation, we get $P_2\Phi(P_1)P_2=0$.
Similarly,  multiplying by $P_1$ from the both sides of this
equation, we have $P_1\Phi(P_2)P_1=0$. For every
$M_{12}\in{\mathcal M}_{12}$,
$\Phi(M_{12})=\Phi(P_1M_{12})=\Phi(P_1)M_{12}+P_1\Phi(M_{12})$ and
so $P_1\Phi(P_1)M_{12}=0$. Hence $P_1\Phi(P_1)P_1=0$. Similarly,
from $\Phi(M_{12})=\Phi(M_{12}P_2)$, one can obtain
$P_2\Phi(P_2)P_2=0$.

Denote $T_0=P_1\Phi(P_1)P_2-P_2\Phi(P_1)P_1$. Define
$\Psi:{\mathcal M}\rightarrow {\mathcal M}$ by
$\Psi(T)=\Phi(T)-[T,T_0]$ for every $T\in {\mathcal M}$. Then it
is easy to see that $\Psi$ is also a nonlinear derivation and
$\Psi(P_1)=\Psi(P_2)=0$. Note that for every $T\in {\mathcal M}:
T\mapsto [T,T_0]$ is an additive derivation of ${\mathcal M}$.
Therefore, without loss of generality, we may assume
$\Phi(P_1)=\Phi(P_2)=0$. Then for every $T_{ij}\in {\mathcal
M}_{ij}$, $\Phi(T_{ij}) =P_i\Phi(M_{ij})P_j\in {\mathcal M}_{ij}$
($i,j=1,2$).

Let $T$ be in $\mathcal M$, write $T=T_{11}+T_{12}+T_{21}+T_{22}$.
In order to prove the additivity of $\Phi$, we only need to show
$\Phi$ is additive on ${\mathcal M}_{ij} (1\leq i,j\leq 2)$ and
$\Phi(T_{11}+T_{12}+T_{21}+T_{22})=\Phi(T_{11})+\Phi(T_{12})+\Phi(T_{21})+\Phi(T_{22})$.
 We will complete the proof by checking two claims.

\textbf{Claim 1.} $\Phi$ is additive on ${\mathcal M}_{ij} (1\leq
i,j\leq 2)$.

Set $T_{ij}, S_{ij}, M_{ij}\in {\mathcal M}_{ij}$.
From $(T_{11}+T_{12})M_{12}=T_{11}M_{12}$, it follows that
$$\Phi(T_{11}+T_{12})M_{12}+(T_{11}+T_{12})\Phi(M_{12})=\Phi(T_{11})M_{12}+T_{11}\Phi(M_{12}).$$ Note that $\Phi(M_{12})\in {\mathcal
M}_{12}$, so $(\Phi(T_{11}+T_{12})-\Phi(T_{11}))M_{12}=0.$ Then
$(\Phi(T_{11}+T_{12})-\Phi(T_{11}))P_1=0$. This implies
$$(\Phi(T_{11}+T_{12})-\Phi(T_{11})-\Phi(T_{12}))P_1=0.$$ Similarly,
from $(T_{11}+T_{12})M_{21}=T_{12}M_{21}$, we have
$(\Phi(T_{11}+T_{12})-\Phi(T_{11})-\Phi(T_{12}))M_{21}=0.$ Then
$$(\Phi(T_{11}+T_{12})-\Phi(T_{11})-\Phi(T_{12}))P_2=0.$$ Thus
$$\Phi(T_{11}+T_{12})=\Phi(T_{11})+\Phi(T_{12}).$$ Similarly, $\Phi(T_{12}+T_{22})=\Phi(T_{12})+\Phi(T_{22}).$
Since
$T_{12}+S_{12}=(P_1+T_{12})(P_2+S_{12})$,
we have that
$$\begin{aligned}
 \Phi(T_{12}+S_{12} )&=\Phi(P_1+T_{12})(P_2+S_{12})+(P_1+T_{12})\Phi(P_2+S_{12})  \\
&=    \Phi(T_{12})+\Phi(S_{12} ). \end{aligned} $$ In the same
way, one can show that $\Phi(T_{21}+S_{21}
)=\Phi(T_{21})+\Phi(S_{21} ).$ That is, $\Phi$ is additive on
${\mathcal M}_{12}, {\mathcal M}_{21}$.

From $(T_{11}+S_{11})M_{12}=T_{11}M_{12}+S_{11}M_{12}$, it follows that
$$\begin{array}{rl}
  & \Phi(T_{11}+S_{11})M_{12}+(T_{11}+S_{11})\Phi(M_{12})\\
= & \Phi(T_{11}M_{12})+\Phi(S_{11}M_{12})\\
= & \Phi(T_{11})M_{12}+
T_{11}\Phi(M_{12})+\Phi(S_{11})M_{12}+S_{11}\Phi(M_{12}).\end{array}$$
Thus  $(\Phi(T_{11}+S_{11})-\Phi(T_{11})-\Phi(S_{11}))M_{12}=0$.
This yields
$$(\Phi(T_{11}+S_{11})-\Phi(T_{11})-\Phi(S_{11}))P_1=0.$$ Note that
$\Phi(T_{11}+S_{11})-\Phi(T_{11})-\Phi(S_{11})\in {\mathcal
M}_{11}$. So $$\Phi(T_{11}+S_{11})=\Phi (S_{11})+\Phi(T_{11}).$$
Similarly ,$\Phi(T_{22}+S_{22})=\Phi (T_{22})+\Phi(S_{22}).$ That
is, $\Phi$ is additive on ${\mathcal M}_{11}, {\mathcal M}_{22}$,
as desired.

\textbf{Claim 2.} $\Phi(T_{11}+T_{12}+T_{21}+T_{22})=\Phi(T_{11})+\Phi(T_{12})+\Phi(T_{21})+\Phi(T_{22})$

From $(T_{11}+T_{12}+T_{21}+T_{22})M_{12}=(T_{11}+T_{21})M_{12}$,
we have
$$\begin{array}{rl}
  & \Phi(T_{11}+T_{12}+T_{21}+T_{22})M_{12}+(T_{11}+T_{12}+T_{21}+T_{22})\Phi(M_{12})\\
= & \Phi(T_{11}M_{12})+\Phi(T_{21}M_{12})\\
= & \Phi(T_{11})M_{12}+ T_{11}\Phi(M_{12})+\Phi(T_{21})M_{12}+T_{21}\Phi(M_{12}).\end{array}$$
Then
$$\begin{array}{ll}
 &(\Phi(T_{11}+T_{12}+T_{21}+T_{22})-\Phi(T_{11})-\Phi(T_{12})-\Phi(T_{21})-\Phi(T_{22}))M_{12}\\
= &
(\Phi(T_{11}+T_{12}+T_{21}+T_{22})-\Phi(T_{11})-\Phi(T_{21}))M_{12}=0.\end{array}$$
This gives
$$(\Phi(T_{11}+T_{12}+T_{21}+T_{22})-\Phi(T_{11})-\Phi(T_{12})-\Phi(T_{21})-\Phi(T_{22}))P_1=0.$$
From $(T_{11}+T_{12}+T_{21}+T_{22})M_{21}=(T_{12}+T_{22})M_{21}$,
it follows that
$$(\Phi(T_{11}+T_{22}+T_{12}+T_{21})-\Phi(T_{11})-\Phi(T_{12})-\Phi(T_{21})-\Phi(T_{22}))P_2=0.$$
So
$\Phi(T_{11}+T_{12}+T_{21}+T_{22})=\Phi(T_{11})+\Phi(T_{12})+\Phi(T_{21})+\Phi(T_{22})$.

Now, we turn to prove Theorem 1.2.

\textbf{Proof of Theorem 1.2.} We will finish the proof of the
Theorem 1.2 by checking several claims.

\textbf{Claim 1.}  $\Phi(0)=0$ and there is $T_0\in {\mathcal M}$ such that $\Phi(P_i)=[P_i,T_0]$ ($i=1,2$).

It is clear that $\Phi(0)=\Phi([0,0]_{\xi})=[\Phi(0),0]_{\xi}+[0,\Phi(0)]_{\xi}=0$.

For every $M_{12}$,
$$\begin{array}{ll}
\Phi(M_{12}) & =\Phi([P_1,M_{12}]_{\xi})=[\Phi(P_1), M_{12}]_{\xi}+[P_1, \Phi(M_{12})]_{\xi}\\
             &  =\Phi(P_1)M_{12}-\xi M_{12}\Phi(P_1)+P_1\Phi(M_{12})-\xi \Phi(M_{12})P_1. \end{array}\eqno{(2)}$$
Multiplying by $P_1,P_2$ from the left and the right in equation
(2) respectively, we have
$$P_1\Phi(P_1)P_1M_{12}=\xi M_{12}P_2\Phi(P_1)P_2. $$ That is
$[P_1\Phi(P_1)P_1+P_2\Phi(P_1)P_2, M_{12}]_{\xi}=0$.  Now Lemma
2.2 yields that $P_1\Phi(P_1)P_1+P_2\Phi(P_1)P_2\in  (\xi
P_1+P_2){\mathcal Z}_{\mathcal M}$. For every $M_{21}$,
$$\begin{array}{ll}
\Phi(M_{21}) & =\Phi([P_2,M_{21}]_{\xi})=[\Phi(P_2), M_{21}]_{\xi}+[P_2, \Phi(M_{21})]_{\xi}\\
             &  =\Phi(P_2)M_{21}-\xi M_{21}\Phi(P_2)+P_2\Phi(M_{21})-\xi \Phi(M_{21})P_2. \end{array}\eqno{(3)}$$
Multiplying by $P_2,P_1$  from the left and the right in equation
(3) respectively, we obtain $$P_2\Phi(P_2)P_2M_{21}=\xi
M_{21}P_1\Phi(P_2)P_1. $$ That is
$[P_2\Phi(P_2)P_2+P_1\Phi(P_2)P_1, M_{21}]_{\xi}=0$. Using Lemma
2.2 again, we get $P_2\Phi(P_2)P_2+P_1\Phi(P_2)P_1\in (P_1+\xi
P_2){\mathcal Z}_{\mathcal M}$. Assume
$P_1\Phi(P_1)P_1+P_2\Phi(P_1)P_2=(\xi P_1+P_2)Z_1$ and
$P_2\Phi(P_2)P_2+P_1\Phi(P_2)P_1=(P_1+\xi P_2)Z_2$, $Z_1,Z_2\in
{\mathcal Z}_{\mathcal M}$.  From $[P_1,P_2]_{\xi}=0$, it follows
that
$$\begin{array}{ll}
  & \Phi([P_1,P_2]_{\xi})=[\Phi(P_1),P_2]_{\xi}+[P_1,\Phi(P_2)]_{\xi}\\
=  & \Phi(P_1)P_2-\xi P_2\Phi(P_1)+P_1\Phi(P_2)-\xi\Phi(P_2)P_1\\
= & (1-\xi)P_1\Phi(P_2)P_1+(1-\xi)P_2\Phi(P_1)P_2+P_1\Phi(P_1)P_2\\
  & +P_1\Phi(P_2)P_2-\xi P_2\Phi(P_2)P_1-\xi P_2\Phi(P_1)P_1\\
= & 0. \end{array}$$ Then
$$P_1\Phi(P_2)P_1=P_2\Phi(P_1)P_2=P_1\Phi(P_1)P_2+P_1\Phi(P_2)P_2=P_2\Phi(P_1)P_1+
P_2\Phi(P_2)P_1=0.\eqno{(4)}$$ A direct computation shows that
$[(\xi
P_1+P_2)Z_1,P_2]_{\xi}=[P_1\Phi(P_1)P_1+P_2\Phi(P_1)P_2,P_2]_{\xi}=0$.
And so $(1-\xi) P_2Z_1=0$.  Then $Z_1MP_2=0$ for all
$M\in{\mathcal M}$. Noting that ${\overline P_2}=I$, we have
$Z_1=0$. That is $P_1\Phi(P_1)P_1+P_2\Phi(P_1)P_2=0$. Similarly,
$P_2\Phi(P_2)P_2+P_1\Phi(P_2)P_1=0$. By (4),
$\Phi(P_1)+\Phi(P_2)=0$. Denote
$T_0=P_1\Phi(P_1)P_2-P_2\Phi(P_1)P_1$. Then it is easy to check
that $T_0$ is the desired.

\vskip 0.4cm{\it Obviously, $T\mapsto [T,T_0]$ is an additive derivation. Without loss of generality,  we may assume that $\Phi(P_1)=\Phi(P_2)=0$. }\vskip 0.4cm

If $\Phi$ is additive, then $\Phi(I)=\Phi(P_1)+\Phi(P_2)=0$.
$\Phi((1-\xi )T)=\Phi([I,T]_{\xi})= [I,\Phi(T)]_{\xi}=(1-\xi
)\Phi(T)$ for all $T\in {\mathcal M}$. So $ \Phi(\xi T)=\xi
\Phi(T)$ for all $T\in {\mathcal M}$. Taking $T,S\in {\mathcal M}$
and noting that $(1-\xi)[S,T]_{-1}=[S,T]_{\xi}+[T,S]_{\xi}$, we
obtain that
$$\begin{array}{ll}
 &\Phi((1-\xi)[S,T]_{-1})=\Phi([S,T]_{\xi})+\Phi([T,S]_{\xi})\\
= & \Phi(S)T-\xi T\Phi(S)+S\Phi(T)-\xi \Phi(T)S+\Phi(T)S-\xi S\Phi(T)+T\Phi(S)-\xi \Phi(S)T\\
= & (1-\xi)( \Phi(S)T+S\Phi(T)+\Phi(T)S+T\Phi(S)).\end{array}$$
Note that $\Phi((1-\xi )T)=(1-\xi )\Phi(T)$ for all $T\in{\mathcal
M}$, it follows that $$\Phi([S,T]_{-1})= [\Phi(S),T]_{-1}+
[S,\Phi(T)]_{-1}$$ for all $T,S\in{\mathcal M}$. Hence $\Phi$ is
an additive Jordan derivation. By \cite{Br3}, $\Phi$ is an
additive derivation which is the conclusion of our Theorem 1.2.
Now we only need to show $\Phi$ is additive. For every $T\in
\mathcal M$, it has the form $T=T_{11}+T_{12}+T_{21}+T_{22}$. Just
like the proof of Theorem 1.1, we will show $\Phi$ is additive on
${\mathcal M}_{ij} (1\leq i,j\leq 2)$ and
$\Phi(T_{11}+T_{12}+T_{21}+T_{22})=\Phi(T_{11})+\Phi(T_{12})+\Phi(T_{21})+\Phi(T_{22})$.
We divide the proof into several steps.

\textbf{Claim 2.} $\Phi(M_{ij})\in {\mathcal M}_{ij}$ for every
$M_{ij}\in {\mathcal M}_{ij}$ ($1\leq i\neq j\leq 2$).

We only treat the case $i=1,j=2$. The other case can be treated similarly.
 Noting $[P_1,M_{12}]_{\xi}=M_{12}$, we have
 $$\begin{aligned}
    \Phi(M_{12})&=\Phi([P_1,M_{12}]_{\xi})
=  [\Phi(P_1), M_{12}]_{\xi}+[P_1, \Phi(M_{12})]_{\xi}\\
&=  [P_1, \Phi(M_{12})]_{\xi}=P_1\Phi(M_{12})-\xi \Phi(M_{12})P_1.
\end{aligned}$$
Then
$$P_2\Phi(M_{12})P_2=P_1\Phi(M_{12})P_1=0. \eqno{(5)}$$ Furthermore, $P_2\Phi(M_{12})P_1=0$, if $\xi\neq -1$, i.e., $\Phi(M_{12})\in {\mathcal M}_{12}$.

Next we treat the case $\xi=-1$. For every $M_{11}$,
$$\begin{aligned}
 \Phi(M_{11}M_{12})&=\Phi([M_{11}, M_{12}]_{-1})
=  [\Phi(M_{11}), M_{12}]_{-1}+[M_{11}, \Phi(M_{12})]_{-1}\\
& = \Phi(M_{11})
M_{12}+M_{12}\Phi(M_{11})+M_{11}\Phi(M_{12})+\Phi(M_{12})M_{11}.\end{aligned}$$
By (5), we have $$P_2\Phi(M_{11}M_{12})P_1=\Phi(M_{12})M_{11}.$$
Then for every $N_{11}$,
$P_2\Phi(N_{11}M_{11}M_{12})P_1=\Phi(M_{12})N_{11}M_{11}$. On the
other hand,
$$P_2\Phi(N_{11}M_{11}M_{12})P_1=\Phi(M_{11}M_{12})N_{11}=\Phi(M_{12})M_{11}N_{11}.$$
Thus $\Phi(M_{12})[N_{11},M_{11}]=0$. For every $R_{11}$,
$$\Phi(M_{12})R_{11}[N_{11},M_{11}]=P_2\Phi(R_{11}M_{12})P_1[N_{11},M_{11}]=0.$$
By Lemma 2.1(ii), $\Phi(M_{12})P_1=0$ which finishes the proof.

\textbf{Claim 3.}  $\Phi(M_{ii})\in {\mathcal M}_{ii}$ for every
$M_{ii}\in{\mathcal M}_{ii}$ $(i=1,2)$.

\textbf{Proof.} Without loss of generality, we only treat the case $i=1$.
$$\begin{array}{ll}
\Phi(P_1) & =\Phi([I,\frac{1}{1-\xi}P_1]_{\xi})=[\Phi(I),\frac{1}{1-\xi}P_1]_{\xi}+[I,\Phi(\frac{1}{1-\xi}P_1)]_{\xi}\\
& =\frac{1}{1-\xi}\Phi([I,P_1]_{\xi})+[I,\Phi(\frac{1}{1-\xi}P_1)]_{\xi}\\
&=\frac{1}{1-\xi}\Phi((1-\xi)P_1)+(1-\xi)\Phi(\frac{1}{1-\xi}P_1)=0.\end{array}$$
Note that $\Phi((1-\xi)P_1)=\Phi([P_1,P_1]_{\xi})=0$, so $\Phi(\frac{1}{1-\xi}P_1)=0$.
$$\begin{aligned}\Phi(M_{11})&=\Phi([\frac{1}{1-\xi}P_1,M_{11}]_{\xi})=[\frac{1}{1-\xi}P_1,\Phi(M_{11})]_{\xi}\\
&=\frac{1}{1-\xi}(P_1\Phi(M_{11})-\xi
\Phi(M_{11})P_1).\end{aligned}$$ This implies $\Phi(M_{11})\in
{\mathcal M}_{11}$.

\textbf{Claim 4.}  For every $T_{ii}$, $T_{ji}$  and $T_{ij}$ $(1\leq i\neq
j\leq 2)$,
$\Phi(T_{ii}+T_{ij})=\Phi(T_{ii})+\Phi(T_{ij}) $,
$\Phi(T_{ii}+T_{ji})=\Phi(T_{ii})+\Phi(T_{ji})$.

Assume $i=1,j=2$. For every $M_{12}\in {\mathcal M}_{12}$, $[T_{11}+T_{12},M_{12}]_{\xi}=[T_{11},M_{12}]_{\xi}$, by Claim 2,
$$[\Phi(T_{11}+T_{12}),M_{12}]_{\xi}+[T_{11}+T_{12},\Phi(M_{12})]_{\xi}=
[\Phi(T_{11}),M_{12}]_{\xi}+[T_{11},\Phi(M_{12})]_{\xi},$$
$$[\Phi(T_{11}+T_{12})-\Phi(T_{11}),M_{12}]_{\xi}=0.$$ From Lemma 2.2,
$$\Phi(T_{11}+T_{12})-\Phi(T_{11})=P_1(\Phi(T_{11}+T_{12})-\Phi(T_{11}))P_2+(\xi
P_1+P_2)Z$$ for some central element $Z\in Z_{\mathcal M}.$
 By computing, $$\begin{aligned}\Phi(T_{12})&=\Phi([P_1,[T_{11}+T_{12},P_2]_{\xi}]_{\xi})\\
 &=[P_1,[\Phi(T_{11}+T_{12}),P_2]_{\xi}]_{\xi}\\
 &=P_1\Phi(T_{11}+T_{22})P_2+{\xi}^2P_2\Phi(T_{11}+T_{22}).\end{aligned}$$
 From Claim 2 and Claim 3, we know that $\Phi(T_{12})=P_1\Phi(T_{11}+T_{12})P_2$ and $P_1\Phi(T_{11})P_2=0$.
 Thus $$\Phi(T_{11}+T_{12})-\Phi(T_{11})=\Phi(T_{12})+(\xi P_1+P_2)Z.$$
Note that
$$\begin{aligned}\Phi([T_{11}+T_{12},P_2]_{\xi})&=[\Phi(T_{11}+T_{12}),P_2]_{\xi}\\
&=[\Phi(T_{11})+\Phi(T_{12})+(\xi
P_1+P_2)Z,P_2]_{\xi}.\end{aligned}$$ On the other hand,
$\Phi([T_{11}+T_{12},P_2]_{\xi})=\Phi([T_{12},P_2]_{\xi})=[\Phi(T_{12}),P_2]_{\xi}$.
Combining this with Claim 3, we have $[(\xi
P_1+P_2)Z,P_2]_{\xi}=0$ and so $ZP_2=0$ which implies $Z=0$.
Similarly, $\Phi(T_{11}+T_{21})=\Phi(T_{11})+\Phi(T_{21})$. The
rest goes similarly.

\textbf{Claim 5.} $\Phi$ is additive on ${\mathcal M}_{12}$ and
${\mathcal M}_{21}$.

Let $T_{12},S_{12}\in{\mathcal M}_{12}$. Since
$T_{12}+S_{12}=[P_1+T_{12}, P_2+S_{12} ]_{\xi}$,
we have that
$$\begin{array}{rl}
 & \Phi(T_{12}+S_{12} )=
     [\Phi(P_1+T_{12}), P_2+S_{12}]_{\xi}+[P_1+T_{12}, \Phi(P_2+S_{12})]_{\xi}  \\
= &  [\Phi(P_1) +\Phi(T_{12}), P_2+S_{12}]_{\xi}+[P_1+T_{12}, \Phi(P_2)+\Phi(S_{12})]_{\xi} \\
= &   \Phi(T_{12})+\Phi(S_{12} ). \end{array} $$
Similarly, $\Phi$ is additive on ${\mathcal
M}_{21}$.

\textbf{Claim 6.} For every $T_{11}\in{\mathcal M}_{11}$,
$T_{22}\in{\mathcal M}_{22}$, $\Phi(T_{11}+T_{22})=\Phi
(T_{11})+\Phi(T_{22})$.

For every $M_{12}\in {\mathcal M}_{12}$, $[T_{11}+T_{22},M_{12}]_{\xi}=T_{11}M_{12}-\xi M_{12}T_{22}$. From Claim 5, it follows that $$\begin{array}{rl}
  & [\Phi(T_{11}+T_{22}),M_{12}]_{\xi}+[T_{11}+T_{22},\Phi(M_{12})]_{\xi}=\Phi([T_{11}+T_{22},M_{12}]_{\xi})\\
= & \Phi(T_{11}M_{12})+\Phi(-\xi M_{12}T_{22})=\Phi([T_{11},M_{12}]_{\xi})+\Phi([T_{22}, M_{12}]_{\xi})\\
= & [\Phi(T_{11}),M_{12}]_{\xi}+
[T_{11},\Phi(M_{12})]_{\xi}+[\Phi(T_{22}), M_{12}]_{\xi}+[T_{22},
\Phi(M_{12})]_{\xi}.\end{array}$$ Thus
$[\Phi(T_{11}+T_{22})-\Phi(T_{11})-\Phi(T_{22}),M_{12}]_{\xi}=0$.
By Lemma 2.2, $$\Phi(T_{11}+T_{22})-\Phi(T_{11})-\Phi(T_{22})\in
{\mathcal M}_{12}+(\xi P_1+P_2){\mathcal Z}_{\mathcal M}.$$ On the
other hand, $[T_{11}+T_{22}, \frac{P_1}{1-\xi}]_{\xi}=T_{11}$.
From the proof of Claim 3, one can see
$\Phi(\frac{P_1}{1-\xi})=0$. Hence $[\Phi(T_{11}+T_{22}),
\frac{P_1}{1-\xi}]_{\xi}=\Phi(T_{11})$, i.e.,
$(1-\xi)\Phi(T_{11})=\Phi(T_{11}+T_{22})P_1-\xi
P_1\Phi(T_{11}+T_{22})$. Multiplying by $P_1$ and $P_2$ from the
left and the right in the above equation, we have
$P_1\Phi(T_{11}+T_{22})P_2=0$. So
$$\Phi(T_{11}+T_{22})-\Phi(T_{11})-\Phi(T_{22})=(\xi P_1+P_2)Z$$ for
some central element $Z\in Z_{\mathcal M}$. Combining
$\Phi(P_1)=0$ and Claim 3, we conclude
$$\begin{array}{ll}
&\Phi([T_{11},P_1]_{\xi})=\Phi([T_{11}+T_{22},P_1]_{\xi})\\
=& [\Phi(T_{11}+T_{22}),P_1]_{\xi}+[T_{11}+T_{22},\Phi(P_1)]_{\xi}\\
= & [\Phi(T_{11})+(\xi P_1+P_2)Z,P_1]_{\xi}.\end{array}$$ Thus
$[(\xi P_1+P_2)Z,P_1]_{\xi}=0$ which implies $Z=0$. This gives
$\Phi(T_{11}+T_{22})=\Phi(T_{11})+\Phi(T_{22})$.

\textbf{Claim 7.} For every $T_{ii},S_{ii}\in {\mathcal M}_{ii}$
$(i=1,2)$, $\Phi(T_{ii}+S_{ii})=\Phi (T_{ii})+\Phi(S_{ii})$.

Assume $i=1$. For every $M_{12}\in {\mathcal M}_{12}$, $[T_{11}+S_{11},M_{12}]_{\xi}=T_{11}M_{12}+S_{11}M_{12}$. From Claim 5, it follows that
$$\begin{array}{rl}
  & [\Phi(T_{11}+S_{11}),M_{12}]_{\xi}+[T_{11}+S_{11},\Phi(M_{12})]_{\xi}=\Phi([T_{11}+S_{11},M_{12}]_{\xi})\\
= & \Phi(T_{11}M_{12})+\Phi(S_{11}M_{12})=\Phi([T_{11},M_{12}]_{\xi})+\Phi([S_{11}, M_{12}]_{\xi})\\
= & [\Phi(T_{11}),M_{12}]_{\xi}+
[T_{11},\Phi(M_{12})]_{\xi}+[\Phi(S_{11}), M_{12}]_{\xi}+[S_{11},
\Phi(M_{12})]_{\xi}.\end{array}$$ Thus
$[\Phi(T_{11}+S_{11})-\Phi(T_{11})-\Phi(S_{11}),M_{12}]_{\xi}=0$.
By Lemma 2.2, $$\Phi(T_{11}+S_{11})-\Phi(T_{11})-\Phi(S_{11})\in
{\mathcal M}_{12}+(\xi P_1+P_2){\mathcal Z}.$$ On the other hand,
Claim 3 tells us that
$P_1(\Phi(T_{11}+S_{11})-\Phi(T_{11})-\Phi(S_{11}))P_2=0$. So
$$\Phi(T_{11}+S_{11})-\Phi (T_{11})-\Phi(S_{11})= (\xi P_1+P_2)Z$$
for some $Z\in {\mathcal Z}_{\mathcal M}$. This further indicates
$$\begin{aligned}
0&=\Phi([T_{11}+S_{11}, P_2]_{\xi})=[\Phi(T_{11}+S_{11}), P_2]_{\xi}\\
&=[\Phi(T_{11})+\Phi(S_{11})+(\xi P_1+P_2)Z, P_2]_{\xi}\\
&=[(\xi P_1+P_2)Z, P_2]_{\xi}.
\end{aligned}$$ Then $P_2Z=0$, consequently, $Z=0$. That is, $\Phi$ is additive
on  ${\mathcal M}_{11}$. Similarly, $\Phi$ is additive on
${\mathcal M}_{22}$.

\textbf{Claim 8.} For every $T_{ii}, T_{jj}, T_{ij}$, $(1\leq
i\neq j\leq 2)$
$\Phi(T_{ii}+T_{jj}+T_{ij})=\Phi(T_{ii})+\Phi(T_{jj})+\Phi(T_{ij})$
.

Assume $i=1, j=2$. For every $M_{12}\in {\mathcal M}_{12}$,
$[T_{11}+T_{22}+T_{12},M_{12}]_{\xi}=[T_{11}+T_{22},M_{12}]_{\xi}$.
By Claim 6, it follows that
$$[\Phi(T_{11}+T_{22}+T_{12}),M_{12}]_{\xi}+[T_{11}+T_{22}+T_{12},\Phi(M_{12})]_{\xi}
=[\Phi(T_{11})+\Phi(T_{22}),M_{12}]_{\xi}
+[T_{11}+T_{22},\Phi(M_{12})]_{\xi}.$$ Thus
$[\Phi(T_{11}+T_{22}+T_{12})-\Phi(T_{11})-\Phi(T_{22}),M_{12}]_{\xi}=0.$
From Lemma 2.2 and Claim 3, we obtain
$$\begin{array}{ll}
&\Phi(T_{11}+T_{22}+T_{12})-\Phi(T_{11})-\Phi(T_{22})\\
= & P_1(\Phi(T_{11}+T_{22}+T_{12})-\Phi(T_{11})-\Phi(T_{22}))P_2+(\xi P_1+P_2)Z\\
= & P_1\Phi(T_{11}+T_{22}+T_{12})P_2+(\xi P_1+P_2)Z\end{array}$$
for some central element $Z$. A direct computation shows that
$$\begin{aligned}
\Phi(T_{12})&=\Phi(P_1(T_{11}+T_{22}+T_{12})P_2)\\
&=\Phi([P_1,[T_{11}+T_{22}+T_{12},P_2]_{\xi}]_{\xi})\\
&=[P_1,[\Phi([T_{11}+T_{22}+T_{12},P_2]_{\xi})]_{\xi}\\
&=P_1\Phi(T_{11}+T_{22}+T_{12})P_2.\end{aligned}$$ Thus
$$\Phi(T_{11}+T_{22}+T_{12})=\Phi(T_{11})+\Phi(T_{22})+\Phi(T_{12})+(\xi
P_1+P_2)Z.$$ It is easy to see
$$\begin{array}{ll}
 & [\Phi(T_{11}+T_{22}+T_{12}),P_2]_{\xi}=\Phi([T_{11}+T_{22}+T_{12},P_2]_{\xi})\\
&=\Phi([T_{12}+T_{22},P_2]_{\xi})=[\Phi(T_{12})+\Phi(T_{22}),P_2]_{\xi}\end{array}.$$
Then $[(\xi P_1+P_2)Z,P_2]=0$, $ZP_2=0$ which implies $Z=0$. That
is
$\Phi(T_{11}+T_{22}+T_{12})=\Phi(T_{11})+\Phi(T_{22})+\Phi(T_{12})$.
The rest goes similarly.

\textbf{Claim 9.} For every $T_{11}, T_{12},T_{21},T_{22}$,
$$\Phi(T_{11}+T_{12}+T_{21}+T_{22})-\Phi(T_{11})-\Phi(T_{12})-\Phi(T_{21})-\Phi(T_{22})\in(\xi
P_1+P_2){\mathcal Z}_{\mathcal M} \cap  (P_1+\xi P_2){\mathcal
Z}_{\mathcal M}.$$ Consequently, if $\xi\neq -1$,
$\Phi(T_{11}+T_{12}+T_{21}+T_{22})=\Phi(T_{11})+\Phi(T_{12})+\Phi(T_{21})+\Phi(T_{22})$.

For every $M_{12}\in {\mathcal M}_{12}$,
$[T_{11}+T_{12}+T_{21}+T_{22},M_{12}]_{\xi}=[T_{11}+T_{21}+T_{22},M_{12}]_{\xi}$.
From Claim 8, it follows that $$\begin{array}{ll}
& [\Phi(T_{11}+T_{12}+T_{21}+T_{22}),M_{12}]_{\xi}+[T_{11}+T_{12}+T_{21}+T_{22},\Phi(M_{12})]_{\xi}\\
= & [\Phi(T_{11})+\Phi(T_{21})+\Phi(T_{22}),M_{12}]_{\xi}
   +[T_{11}+T_{21}+T_{22},\Phi(M_{12})]_{\xi}.\end{array}$$ Thus
$$[\Phi(T_{11}+T_{12}+T_{21}+T_{22})-\Phi(T_{11})-\Phi(T_{21})-\Phi(T_{22}),M_{12}]_{\xi}=0.$$
Since $\Phi(T_{12})\in {\mathcal M}_{12}$, we have
$$[\Phi(T_{11}+T_{12}+T_{21}+T_{22})-\Phi(T_{11})-\Phi(T_{12})-\Phi(T_{21})-\Phi(T_{22}),M_{12}]_{\xi}=0.$$
Similarly, from
$[T_{11}+T_{12}+T_{21}+T_{22},M_{21}]_{\xi}=[T_{11}+T_{12}+T_{22},M_{21}]_{\xi},$
we can obtain
 $$[\Phi(T_{11}+T_{12}+T_{21}+T_{22})-\Phi(T_{11})-\Phi(T_{12})-\Phi(T_{21})-\Phi(T_{22}),M_{21}]_{\xi}=0.$$
 From Lemma 2.2, it follows that
  $$\Phi(T_{11}+T_{12}+T_{21}+T_{22})-\Phi(T_{11})-\Phi(T_{12})-\Phi(T_{21})-\Phi(T_{22})\in(\xi P_1+P_2){\mathcal Z}_{\mathcal M} \cap  (P_1+
  \xi P_2){\mathcal Z}_{\mathcal M}.$$

Note that if $\xi\neq -1$,   $(\xi P_1+P_2){\mathcal Z}_{\mathcal M} \cap  (P_1+\xi P_2){\mathcal Z}_{\mathcal M}=\{0\}$. Thus
$$\Phi(T_{11}+T_{12}+T_{21}+T_{22})=\Phi(T_{11})+\Phi(T_{12})+\Phi(T_{21})+
\Phi(T_{22}).$$

\textbf{Claim 10.} If $\xi=-1$, $\Phi(T_{11}+T_{12}+T_{21}+T_{22})=\Phi(T_{11})+\Phi(T_{12})+\Phi(T_{21})+
\Phi(T_{22})$ holds true, too.

By Claim 9, we may assume
$\Phi(T_{12}+T_{21})=\Phi(T_{12})+\Phi(T_{21})+(-P_1+P_2)Z_1$, $\Phi(T_{11}+T_{12}+T_{21})=\Phi(T_{11})+\Phi(T_{12})+\Phi(T_{21})+(-P_1+P_2)Z_2$
 and $\Phi(T_{11}+T_{12}+T_{21}+T_{22})=\Phi(T_{11})+\Phi(T_{12})+\Phi(T_{21})+
\Phi(T_{22})+(-P_1+P_2)Z_3$. The following is devoted to showing
$Z_1=Z_2=Z_3=0$. Since
$\Phi(T_{12}+T_{21})=\Phi([T_{12}+T_{21},P_1]_{-1})=[\Phi(T_{12}+T_{21}),P_1]_{-1}$,
substituting
$\Phi(T_{12}+T_{21})=\Phi(T_{12})+\Phi(T_{21})+(-P_1+P_2)Z_1$ into
above equation, we have
 $(-P_1+P_2)Z_1=[(-P_1+P_2)Z_1,P_1]_{-1}=-2P_1Z_1$. Then $Z_1P_1=Z_1P_2=0$ and so $Z_1=0$.
From
$$\begin{array}{ll}
&[\Phi(T_{11})+\Phi(T_{12})+\Phi(T_{21})+(-P_1+P_2)Z_2,P_2]_{-1}\\
= &[\Phi(T_{11}+T_{12}+T_{21}),P_2]_{-1}=\Phi([T_{11}+T_{12}+T_{21},P_2]_{-1})\\
=&\Phi(T_{12}+T_{21})=\Phi(T_{12})+\Phi(T_{21})\\
=&[\Phi(T_{11})+\Phi(T_{12})+\Phi(T_{21}),P_2]_{-1},\end{array}$$
it follows that $[(-P_1+P_2)Z_2, P_2]=0$. Thus $Z_2=0$. At last,
$$\begin{array}{ll}
 & [\Phi(T_{11})+\Phi(T_{12})+\Phi(T_{21})+\Phi(T_{22})+(-P_1+P_2)Z_3,P_1]_{-1} \\
 = & [\Phi(T_{11}+T_{12}+T_{21}+T_{22}),P_1]_{-1}=\Phi([T_{11}+T_{12}+T_{21}+T_{22},P_1]_{-1})\\
=
&\Phi([T_{11}+T_{12}+T_{21},P_1]_{-1})=[\Phi(T_{11})+\Phi(T_{12})+\Phi(T_{21})+\Phi(T_{22}),P_1]_{-1}.\end{array}$$
So $[(-P_1+P_2)Z_3, P_1]_{-1}=0$ which implies $Z_3=0$. Hence
$\Phi(T_{11}+T_{12}+T_{21}+T_{22})=\Phi(T_{11})+\Phi(T_{12})+\Phi(T_{21})+
\Phi(T_{22})$, as desired.


\begin{thebibliography}{99}

\bibitem{AH} R. An, J. Hou, Additivity of Jordan multiplicative maps
on Jordan operator algebras, Taiwanese J. Math., 10(2006), 45-64.

\bibitem{BD} Z.F. Bai, S.P. Du and J.C. Hou, Multiplicative Lie isomorphisms between prime rings, Communications in Algebra, 36(2008), 1626-1633.

\bibitem{Br} M. Bre\v{s}ar, Centralizing mappings on von Neumann algebras, Proc. Ams. Math. Soc., 111(1991), 501-510.

\bibitem{Br2} M. Bre\v{s}ar, Commuting traces of biadditive mappings,
commutativity preserving mappings, and Lie mappings, Trans. Amer.
Math. Soc., 335(1993), 525-546.

\bibitem{Br3} M. Bre\v{s}ar, Jordan derivation on semiprime rings,  Proc. Ams. Math. Soc., 104(1988), 1003-1006.

\bibitem{Br4} M. Bre\v{s}ar, Jordan derivations revised, Math. Proc. Cambridge Philos. Soc., 139(2005), 411-425.

\bibitem{BrBuPe} J.A. Brooke, P. Brusch, B. Pearson, Commutativity up to a factor of bounded operators in complex Hilbert spaces,
Roy. Soc. Lond. Proc. Ser. A Math. Phy. Eng. Sci. A, 458(2002),
109-118.

\bibitem{Ch} E. Christensen, Derivations of nest algebras, Ann. Math., 229(1977), 155-161.

\bibitem{John} B.E. Johnson, Symmetric amenability and the nonexistence of Lie and Jordan derivations, Math. Proc. Cambridge Philos. Soc., 120(1996), 455-473.

\bibitem{Kad} R.V. Kadison and J.R. Ringrose, Fundamentals of the
theory of operator algebras, Vol. I, Acdemic Press, New York,
1983; Vol. II Acdemic Press, New York, 1986.

\bibitem{Ka} C. Kassel. Quantum Group, Springer-verlag, New York, 1995.

\bibitem{MaVi} M. Mathieu, A.R. Villena, The structure of Lie derivations on C*-algebrs, J. Funct. Anal., 202(2003), 504-525.

\bibitem{Lu1} F. Lu, Multiplicative mappings of operator algebras,
Linear Algebra Appl., 347(2002), 283-291

\bibitem{Mar} W.S. Martindale III, When are multiplicative mappings additive?, Proc. Amer. Math. Soc., 21(1969), 695-698.

\bibitem{Mie2} C.R. Miers, Lie isomorphisms of operator algebras, Pacific J. of
Math., 38(1971), 717-735.

\bibitem{Mie3} C.R. Miers, Lie derivations of von Neumann algebras, Duke Math. J., 40(1973), 403-409.

\bibitem{QiHo1} X.F. Qi, J. Hou, Characterizations of $\xi$-Lie multiplicative isomorphisms,
Proceeding of the 3rd International workshop of Matrix analysis and Applications, 2009.

\bibitem{QiHo2} X.F. Qi, J. Hou, Additive Lie ($\xi$-Lie) derivations and generalized Lie ($\xi$-Lie) derivations on nest algebras,
 Linear Algebra Appl., 431(2009), 843-854.

\bibitem{QH} X.F. Qi, J. Hou, Characterization of Lie multiplicative isomorphisms between nest
 algebras, Science China Mathematics, 54(2011), 2453-2462.


\bibitem{Sa} S. Sakai, Derivations of $W^*$-algebras, Ann. Math., 83(1966), 273-279.

\bibitem{Se1} P. \v{S}emrl, Additive derivations of some operator algebras, Illinois J. Math., 35(1991), 234-240.

\bibitem{Se2} P. \v{S}emrl, Rings derivations on standard operator algebras, J. Funct. Anal., 112(1993), 318-324.

\bibitem{Wang} Y. Wang, Additivity of multiplicative maps on
triangular rings, Linear Algebra Appl., 434(2011), 625-635.

\bibitem{YuZh} W.Y. Yu, J.H. Zhang, Nonlinear Lie derivations of triangular algerbas, Linear Algebra Appl., 432(2010), 2953-2960.

\end{thebibliography}
\end{document}